\documentclass{amsart}

\usepackage {amsfonts}
\usepackage[english]{babel}
\def\g{\gamma}
\def\G{\Gamma}
\def\d{\delta}
\def\a{\alpha}
\def\b{\beta}
\def\p{\varphi}
\def\e{\varepsilon}
\def\l{\lambda}
\def\L{\Lambda}
\def\s{\sigma}
\def\k{\kappa}
\def\o{\omega}
\def\O{\Omega}
\def\H{\mathcal H}

\def\R{{\mathbb R}}
\def\C{{\mathbb C}}
\def\N{{\mathbb N}}
\def\Z{{\mathbb Z}}

\def\Im{\mbox{Im }}
\def\bs{~\hfill\rule{7pt}{7pt}}

\DeclareMathOperator{\spec}{sp}
\DeclareMathOperator{\Res}{Res}
\DeclareMathOperator{\const}{const}
\DeclareMathOperator{\supp}{supp }
\DeclareMathOperator{\dist}{dist}

\newtheorem{Th}{Theorem}
\newtheorem{Pro}{Proposition}
\newtheorem{Def}{Definition}
\newtheorem{Cor}{Corollary}

\begin{document}

\title{Non-negative crystalline and Poisson measures in the Euclidean space}

\author{Sergii Yu. Favorov}

\address{Sergii Favorov,
\newline\hphantom{iii}  V.N.Karazin Kharkiv National University
\newline\hphantom{iii} Svobody sq., 4, Kharkiv, Ukraine 61022}
\email{sfavorov@gmail.com}

\maketitle {\small
\begin{quote}
\noindent{\bf Abstract.}
We study  properties of temperate non-negative purely atomic measures in the Euclidean space such that  the distributional Fourier transform of these measures are pure point ones.
A connection between these measures and almost periodicity is shown,  several forms of the uniqueness theorem are proved. We also obtain necessary and sufficient conditions for a measure
with positive integer masses on the real line to correspond the zero set of an absolutely convergent Dirichlet series with bounded spectrum.
\medskip


\end{quote}
}

\medskip

   \section{Introduction}\label{S1}
   \bigskip

Following \cite{M1}, we say that a pure point complex measure $\mu$ on $\R^d$ is a Poisson measure if $\mu$ is a temperate distribution and its Fourier transform $\hat\mu$ is also a pure point measure. If $\mu$ and $\hat\mu$ have locally finite supports, then $\mu$ is said to be a crystalline measure. If, in addition, the measures $|\mu|$ and $|\hat\mu|$ are temperate distributions, then $\mu$ is called the Fourier quasicrystal.
\smallskip

Fourier quasicrystals and crystalline measures are currently being studied very actively. Many works are devoted to the study of the properties of Fourier quasicrystals.
For example, see the collections of papers \cite{D}, \cite{Q}, in particular, the basic work \cite{Lag1}. Fourier quasicrystals have applications in modern physics,
where they are used as mathematical models of certain atomic structures.

Measures of the form
\begin{equation}\label{a}
\mu_A=\sum_n \d_{a_n},\quad A=\{a_n\}\subset\R^d,
\end{equation}
are the most important case of Fourier quasicrystals. A nontrivial example of such measures in the case $d=1$,  where the sequence $A$ is not a finite union of shifts
of an arithmetical progression, was found by P.Kurasov and P.Sarnak \cite{KS}.
Then A.Olevskii and A.Ulanovskii in \cite{OU}, \cite{OU1} showed that a sequence $A\subset\R$, which generate Fourier quasicrystals \eqref{a} is the zero set  of some exponential polynomial
\begin{equation}\label{s}
P(z)=\sum_{1\le j\le N} q_n  e^{2\pi iz\o_j}, \qquad q_j\in\C,\quad \o_j\in\R,
\end{equation}
with real zeros and, conversely, every  zero set $A$  of each exponential polynomial \eqref{s} with only real zeros generates Fourier quasicrystal \eqref{a}.
A multidimensional version of the last result was obtained by W.Lawton and A.Tsikh \cite{Law}.

Here we suppose that $A$ is a multiset, that is, the multiplicity of each point $a_n$ in \eqref{a} (the mass of $\mu_A$ at $a_n$) is finite, and equals the multiplicity of zero $P(z)$ at this point.
\smallskip

However, the condition that the support of $\hat\mu$ is locally finite often too restrictive, so more general objects, Poisson measures, arise (cf. \cite{M2}).
Note that the class of Poisson measures is the widest class inducing the general Poisson formula: if
$$
\mu=\sum_{\l\in\L}a_\l\d_\l
$$
 and
 $$
 \hat\mu=\sum_{\g\in\G}b_\g\d_\g
 $$
are measures on $\R^d$ with countable $\L,\,\G\subset\R^d$, we have for every Schwartz function $\p$
$$
  \sum_{\l\in\L}a_\l\p(\l)=\sum_{\g\in\G}b_\g\hat\p(\g).
$$
Most effects with Fourier quasicrystals are associated with their almost periodicity, so our study of Poisson measures will begin with the following theorem:

\begin{Th}\label{T1}
For every non-negative Poisson measure $\mu$ on $\R^d$ both convolutions $\mu\star\psi$ with any continuous compactly supported function $\psi$ and $\hat\mu\star\p$ with any Schwartz function $\p$
are almost periodic functions in the sense of Bohr.
\end{Th}
Using definitions of the next section, the measure $\mu$ is  almost periodic, and the measure $\hat\mu$ is an almost periodic distribution.

Note that every complex Fourier quasicrystal is an almost periodic distribution \cite[Lemma 1]{F4}. But there are crystalline measures that are not almost periodic distributions \cite{F2}.

\smallskip

 Theorem \ref{T1} implies various properties of non-negative Poisson measures.
\smallskip

Let
\begin{equation}\label{ball}
B(x_n,R_n)=\{x\in\R^n:\,|x-x_n|<R_n\},\quad n=1,2,\dots
\end{equation}
be a sequence of arbitrary balls such that $ R_n\to\infty$ as $n\to\infty$.
\begin{Th}\label{T2}
Let $\mu,\,\nu$ be two non-negative Poisson measures on $\R^d$. If either $|\mu-\nu|B(x_n,R_n)=o(R^d_n)$, or $|\hat\mu-\hat\nu|B(x_n,R_n)=o(R_n^d)$ as $ n\to\infty$, then $\mu\equiv\nu$.
\end{Th}

\begin{Cor}\label{C1}Let $\mu$ be a non-negative Poisson measure on $\R^d$. If $\mu(B(x_n,R_n))=o(R^d_n)$ or $|\hat\mu|(B(x_n,R_n)=o(R_n^d)$ as $ n\to\infty$, then $\mu\equiv0$.
\end{Cor}
In particular, condition of this theorem are fulfilled when $\supp\hat\mu$ does not intersect with an open cone in $\R^d$,
hence any lighthouse $(\mu,S)$ (see \cite{M2}) with a closed cone $S$ could not be a Poisson measure.
\smallskip

Let $E=\cup_n B(x_n,R_n)$ be the union of balls \eqref{ball}, where $R_n\to\infty$ as $n\to\infty$.
\begin{Th}\label{T3}
Let
$$
\mu_A=\sum_{a_n\in A}\d_{a_n},\quad \mu_C=\sum_{c_n\in C}\d_{c_n},\quad A,\ C\subset\R^d,
$$
be two Poisson measures such that
under  appropriate numbering $a_n-c_n\to0$ as $n\to\infty,\,n\in H$, where $H=\{n\in\N:\,a_n\in E\ \text{or}\ c_n\in E\}$.  Then $\mu_A\equiv\mu_C$.
\end{Th}

Other forms of uniqueness theorems can be found in \cite{F1}. They were formulated  for complex Fourier quasicrystals, but  they are actually valid for Poisson measures $\mu$, for which $|\hat\mu|$
are temperate distributions. Indeed, the last condition implies that $\mu$ is an  almost periodic distributions (see \cite[Lemma 1]{F4}), and all proofs are based only on this property.

\begin{Th}\label{T4}
Let $\mu$ be a non-negative Poisson measure on $\R^d$, let $E$ be as above, and assume that
 for all $\l,\,\l'\in E\cap\supp\mu,\,\l\neq\l'$ the conditions $\mu(\l)\ge c>0$ and $|\l-\l'|\ge\b>0$ hold. Assume also that
\begin{equation}\label{mu}
|\hat\mu|(B(0,r))=O(r^d)\qquad (r\to\infty).
\end{equation}
Then for some $\l_j\in\supp\mu$ and nondegenerated linear operators $T_j$ in $\R^d$
 \begin{equation}\label{m}
   \supp\mu=\cup_{j=1}^N T_j\Z^d+\l_j,\qquad\mu=\sum_{j=1}^N\sum_{\l\in T_j\Z^d+\l_j} Q_j(\l)\d_\l,
\end{equation}
where $Q_j$ are Dirichlet series of the form
\begin{equation*}
Q_j(x)=\sum_{\o\in\O_j} q_\o^{(j)}  e^{2\pi ix\cdot\o}, \qquad q_\o^{(j)}\in\C,\quad \sum_{\o\in\O_j} |q_\o^{(j)}|<\infty,\quad \O_j\ \text{are bounded subsets of }\ \R^d.
\end{equation*}
If we replace the condition  $\mu(\l)\ge c>0$ with  $|\mu(\l)|\ge c>0$ for $\l\in E\cap\supp\mu$, then \eqref{m}  is also true for complex Poisson measures.
\end{Th}
Here we set
$$
x\cdot\o=x_1\o_1+\dots x_d\o_d\quad\text{for}\quad x=(x_1,\dots,x_d)\in\R^d,\ \o=(\o_1,\dots,\o_d)\in\R^d.
$$

\medskip
It was shown in \cite{F} that for $d=1$ the zero set $A=\{a_n\}$ of Dirichlet series \eqref{s} under condition $A\subset\R$ generates a Poisson measure $\mu_A=\sum_n\d_{a_n}$.
Here we prove the precise result for dimension $1$:

\begin{Th}\label{T5}
In order that a sequence $A=\{a_n\}\subset\R$, which generate the  measure $\mu_A$ in  \eqref{a}  be the zero set of some  Dirichlet series
\begin{equation}\label{S}
Q(x)=\sum_{\o\in\O} q_\o e^{2\pi ix\o}, \qquad q_\o\in\C,\quad \sum_{\o\in\O_j} |q_\o|<\infty,\quad \O\ \text{is bounded countable subsets of }\ \R,
\end{equation}
 with real zeros it is necessary and sufficient that the Fourier transform $\hat\mu_A$  be a pure point measure
\begin{equation}\label{b}
  \hat\mu_A=\sum_{\g\in\G}b_\g\d_\g
\end{equation}
  such that
\begin{equation}\label{int}
  \int_0^1 t^{-1}|\hat\mu_A|(dt)=\sum_{0<\g<1}\frac{|b_\g|}{\g}<\infty
\end{equation}
and
\begin{equation}\label{var}
	\log|\hat\mu|(B(0,r))= o(r) \qquad (r\to\infty).
\end{equation}
 Here a multiplicity of every zero $a_n$ of $Q(x)$ coincides with the mass of the measure $\mu_A$ at this point.

For the part of sufficiently condition \eqref{var} can be replaced with $\log|\hat\mu|(B(0,r))=O(r)$ .
\end{Th}
Note that a non-negative measure  is a temperate distribution if and only if its masses in the balls $B(0,r)$ grow at most polynomially as $r\to\infty$ (cf.\cite{F1}).
Hence condition \eqref{var} is weaker than $|\hat\mu|$ to be a temperate distribution.

\medskip

The article is structured as follows.

In Section \ref{S2} we give all necessary definitions and formulate auxiliary theorems on almost periodic functions and sets, Fourier transforms, Dirichlet series,  entire functions of exponential growth.

In Section \ref{S3} we show some properties of temperate measures and prove  Theorems \ref{T1}--\ref{T4}.

In Section \ref{S4} we prove the sufficiency part of Theorem \ref{T5}.

In Section \ref{S5} we  prove the necessity part of Theorem \ref{T5}.

\bigskip
\section{Preliminary results}\label{S2}
\bigskip

{\bf The Fourier transform in the sense of distributions}

Denote by $S(\R^d)$ the Schwartz space of test functions $\p\in C^\infty(\R^d)$ with the finite norms
$$
N_n(\p)=\sup_{\R^d}\max_{\|k\|\le n}\,(\max\{1,|x|^n\}) |D^k\p(x)|,\quad n=0,1,2,\dots,
$$
where
$$
|x|=(x_1^2+\dots+x_d^2)^{1/2},\   k=(k_1,\dots,k_d)\in(\N\cup\{0\})^d,\ \|k\|=k_1+\dots+k_d,\  D^k=\partial^{k_1}_{x_1}\dots\partial^{k_d}_{x_d}.
$$
These norms generate the topology on $S(\R^d)$.  Elements of the space $S^*(\R^d)$ of continuous linear functionals on $S(\R^d)$ are called {\it temperate distributions}.
If $\mu\in S^*(\R^d)$
is a measure, we will say that $\mu$ is a {\it temperate measure}.
The Fourier transform of a temperate distribution $f$ is defined by the equality
$$
\hat f(\p)=f(\hat\p)\quad\mbox{for all}\quad\p\in S(\R^d),
$$
where
$$
\hat\p(t)=\int_{\R^d}\p(x)e^{-2\pi i x\cdot t}dx
$$
is the Fourier transform of the function $\p$. By $\check\p$ we  denote the inverse Fourier transform of $\p$.
The Fourier transform is a continuous bijection of $S(\R^d)$ on itself and a bijection of $S^*(\R^d)$ on itself.
\medskip

{\bf Dirichlet series on $\R^d$}

Denote by $\H$ the algebra of all  Dirichlet series
$$
Q(x)=\sum_n r_ne^{2\pi i\o_n\cdot x},\,\o_n\in\R^d,
$$
with  finite Wiener's norm $\|Q\|_W=\sum_n|r_n|$.

For any $Q\in\H$ and any analytic function $h(z)$ on a neighborhood of the set $\overline{\{Q(x):\,x\in\R^d\}}$ we have $h(Q(x))\in\H$ (see \cite{R1}, Ch.VI).
In particular,

for every $Q\in\H$ we have $\exp Q\in\H$,

if $\inf_{\R^d} |Q(x)|>0$, we have $1/Q\in\H$,

if $|Q(x)|\le C<1$, we have $\log(1+Q(x))\in\H$.
\smallskip

\smallskip
{\bf Entire functions  of exponential growth}

Zeros $a_n$ of any entire function $g(z),\,z\in\C$, of exponential growth satisfy the conditions
\begin{equation}\label{l1}
	\#\{a_n:\,|a_n|\le r\}=O(r)\quad\text{as}\quad r\to\infty,
\end{equation}
(as usual, $\# E$ means the number of points of a finite multiset $E$, taking into account a multiplicity of each point), and
\begin{equation}\label{l2}
	\sum_{n:\,|a_n|\le r,a_n\neq0}\frac{1}{a_n}=O(1)\quad\text{as}\quad r\to\infty.
\end{equation}
On the other hand, if a sequence $A=\{a_n\}\subset\C\setminus\{0\}$ satisfies these conditions, then the function
\begin{equation}\label{F}
	F(z)= \prod_{n=-\infty}^\infty(1-z/a_n)e^{z/a_n}
\end{equation}
is an entire function of exponential growth (Lindel\"of's Theorem, see \cite{K}). Note that each $a\in\C$  can be repeated in the sequence $\{a_n\}$ any finite number of times.
\smallskip

If a function $g(z)$ of the exponential type $\s$ is bounded on the real axis, then for all $z\in\C$
$$
|g(z)|\le e^{\s|\Im z|}\sup_{z\in\R}|g(z)|\quad\text{ (Phragment--Lindel\"of Principle, see \cite{K})}.
$$
\medskip

{\bf Almost periodic functions and measures}

\begin{Def}[cf. \cite{B}, \cite{M}]\label{D1}
 A continuous function $g(x)$ on $\R^d$
is  almost periodic if for any  $\e>0$ the set of $\e$-almost periods
  $$
E_{\e}= \{\tau\in\R^d:\,\sup_{x\in\R^d}|g(x+\tau)-g(x)|<\e\}
  $$
is relatively dense, i.e., $E_{\e}\cap B(x,L)\neq\emptyset$ for all $x\in\R^d$ and some $L$ depending on $\e$.
\end{Def}
For example every sum
\begin{equation}\label{Q}
Q(x)=\sum q_n e^{2\pi ix\cdot\o_n},\quad\o_n\in\R^d,\quad q_n\in\C,
\end{equation}
with  finite Wiener's norm $\|Q\|_W=\sum_n |q_n|$ is an almost periodic function.
Note that each almost periodic function is uniformly bounded and uniformly continuous. Also, a uniform in $\R^d$ limit of a sequence of almost periodic functions is almost periodic.

Set $U(T)=\{x=(x_1,\dots,x_d)\in\R^d:\,-T/2<x_j\le T/2,\,j=1,\dots,d\}$. We can prove the existence of the mean value of any almost periodic function
$$
  Mg=\lim_{T\to\infty}\frac{1}{T^d}\int_{U(T)}g(x)dx
$$
the same method as in the dimension $1$ in \cite{B}, also see \cite{FK}.
The Fourier coefficient
$$
  b_\o(g)=\lim_{T\to\infty}\frac{1}{T^d}\int_{U(T)} g(x) e^{-2\pi i\o\cdot x}dx.
$$
 is assigned to every almost periodic function $g$ and every $\o\in\R^d$; spectrum  of  $g$ is the set
 $$
 \spec g=\{\o:\, b_\o(g)\neq0\}.
 $$
It is easy to see that  $b_{\o_n}(Q)=q_n$ and $\spec Q(x)=\{\o_n:\,q_n\neq0\}$. Note that spectra of almost periodic functions are at most countable.
\medskip

\begin{Def}[cf.\cite{R}]\label{D2}
 A measure $\mu$ is almost periodic, if for any continuous function $\p$ with compact support the convolution of $\mu\star\p=\int\p(t-x)\mu(dx)$  is  almost periodic in $t\in\R^d$.
\end{Def}
\begin{Def}[cf.\cite{R}]\label{D3}
 A distribution $F$  is almost periodic, if for any function $\p\in S(\R^d)$ the convolution of $\p\star F=(F(x),\p(t-x))$  is  almost periodic in $t\in\R^d$.
\end{Def}
Note that there exist non-almost periodic measures that are almost periodic in sense of distributions (\cite{M2}, \cite{FK1}).
\bigskip

{\bf Almost periodic sets}

Here we give the definition of  almost periodic sets only for sets on $\R$. The original geometric definition for sets on any horizontal strip  in $\C$ see
\cite[Appendix VI]{L}, \cite{T}. Comparison of various definitions of almost periodic sets see \cite{FRR}.  Almost periodic sets in $\R^d$ were considered in \cite{FK}, \cite{FK1}.

\begin{Def}[\cite{FRR}]\label{D4}
Let $A=\{a_n\}\subset\R$ be a sequence without finite limit points. If the convolution of $\mu_A=\sum_n\d_{a_n}$ with any $C^\infty$-function with compact support is  almost periodic,
 then $A$ is called an almost periodic set.
\end{Def}
It is easy to see that any almost periodic set is relatively dense. Then every point $x\in\R$ can appears in $A$  any finite number of times, hence almost periodic sets are in fact multisets.
\smallskip

A zero set of every Dirichlet series \eqref{Q} with only real zeros is almost periodic \cite{L}.
\footnote{In fact, the zero set of any almost periodic holomorphic function in a strip in $\C$ is almost periodic.
 Moreover, the divisor of any almost periodic holomorphic function in a tube subdomain in $C^d$ is also almost periodic. The corresponding definitions of almost periodic functions and sets in a strip
and connections between these objects see \cite{B}, \cite{L}, \cite{T}, \cite{F0}. Multidimensional case see \cite{R}, \cite{FRR1}, \cite{Fa}.}
\smallskip

 Let $A\subset\R$ be an almost periodic set. Take any $C^\infty$-function $\p\ge0$ with compact support such that $\p(x)\equiv1$ for $0<x<1$. Since $\mu_A\star\p$ is almost periodic, it is uniformly bounded.
 We get
\begin{Pro}\label{P1}
Numbers $\#(A\cap(x,x+1))$  are uniformly bounded in $x\in\R$.
\end{Pro}
 Since any continuous function $g$ with compact support  can be approximated in the uniformly metric with a sequence $\p_n\in C^\infty$ with compact supports such that $\cap_n\supp\p_n=\supp g$, it follows from  Proposition \ref{P1} that  we can replace $C^\infty$-functions by continuous functions with compact support in  Definition \ref{D4}.

Just as the existence of the mean value of an almost periodic function in \cite{B} was proven, using  Proposition \ref{P1} it is easy to show the existence of the density
 $$
D=\lim_{l\to\infty}\frac{\# A\cap(x,\,x+l)}{l}
$$
uniformly with respect to $x\in\R$ (cf.\cite{FK}). Clearly, $D>0$.

It follows from \cite[Lemma 5.5]{M1} that in the case $\hat\mu_A=\sum_{\g\in\G}b_\g\d_\g$  we get  $D=b_0$.
\begin{Pro}[\cite{F}]\label{P2}
 Let $A=\{a_n\}\subset\R$ be an almost periodic set of density $D$ such that $a_n\le a_{n+1}$ for all $n\in\Z$.
Then
$$
a_n=n/D+\phi(n)\quad\text{with an almost periodic function}\quad \phi:\,Z\to\R.
$$
Conversely, any set $A=\{kn+\phi(n)\}$, $\phi$ as above, form an almost periodic set with density $1/k$.
\end{Pro}

\bigskip
\section{Properties of temperate  measures and proofs of theorems \ref{T1}--\ref{T4}}\label{S3}
\bigskip

\begin{Pro}\cite[Lemma 1]{F4}\label{P3}.
Let $\mu$ be a temperate measure and $\hat\mu$ be measure such that $|\hat\mu|\in S^*(\R^d)$.  Then $\mu$ is an almost periodic distribution.
\end{Pro}
 The  condition $|\hat\mu|\in S^*(\R^d)$ is equivalent to the polynomially growth of masses $|\hat\mu|$ in the balls $B(0,r)$ as $r\to\infty$ (cf.\cite{F1}).
\medskip

We will say that a measure $\mu$ is translation bounded, if $\sup_{x\in\R^d}|\mu|B(x,1)<\infty$. It is easy to see that any  translation bounded measure is a temperate measure.
Every almost periodic measure is translation bounded (see \cite{R} or \cite{M}).

\begin{Pro}\cite[Lemma 5.8]{M1}\label{P4}.
If $\mu$ is a translation bounded measure and $\hat\mu$ is a pure point measure, then $\mu$ is almost periodic.
\end{Pro}
This Meyer's proposition is of fundamental importance in our investigation. For convenience of readers, we give its complete proof here.
\smallskip

{\bf Proof}. Let $\hat\mu=\sum_{\g\in\G}b_\g\d_\g$ and $\p$ be any $C^\infty$-function on $\R^d$ with compact support. We have
$$
\mu\star\check\p(t)=(\hat\mu(x),\p(x)e^{2\pi it\cdot x})
=\sum_{\g\in\G\cap\supp\p}b_\g\p(\g) e^{2\pi it\cdot \g}
$$
Since $\hat\mu$ is a measure, we get
$$
\sum_{\g\in\G\cap\supp\p}|b_\g||\p(\g)|\le\sup|\p(\g)|\sum_{\g\in\G\cap\supp\p}|b_\g|<\infty.
$$
Hence the function
$\mu\star\check\p(t)$ is almost periodic.

Now let $\phi$ is an arbitrary function from $S(\R^d)$. Since
$\hat\phi\in S(\R^d)$, we can approximate the function $\check\phi(x)=\hat\phi(-x)$ with compactly supported $C^\infty$-functions $\p_n$ in the metric of the space $S(\R^d)$.
Therefore the functions  $\hat\p_n$ approximate the function $\phi$ in the same metric.  Also,
$$
|\phi(x)-\hat\p_n(x)|\le N_{d+1}(\phi-\hat\p_n)\min\{1,|x|^{-d-1}\}.
$$
Therefore,
$$
 |\mu\star\phi(t)-\mu\star\hat\p_n(t)|\le
 N_{d+1}(\phi-\hat\p_n)\int\min\{1,|x-t|^{-d-1}\}|\mu|(dx)=
$$
$$
 N_{d+1}(\phi-\hat\p_n)|\mu|B(t,1)+ N_{d+1}(\phi-\hat\p_n)\int_1^\infty r^{-d-1}d|\mu|(B(t,r)).
$$
Since the measure $\mu$ is translation bounded, we see that $|\mu|B(t,1)\le C_0$ and $|\mu|B(t,r)\le C_d r^d$ as $r>1$ uniformly with respect to $t\in\R^d$. Therefore the both terms tend to $0$ uniformly with respect to $t$ as
$\hat\p_n\to\phi$ in the space $S(\R^d)$, and the convolution $\phi\star\mu(t)$ is almost periodic. For any continuous function
$\psi$ with compact support there  are functions $\phi_n\in S(\R^d)$ with compact supports, which uniformly converge to $\psi$ and $\cap_n\supp\phi_n=\supp\psi$. Taking into account that $\mu$ is translation bounded,
we get again that the functions $\phi_n\star\mu(t)$ tend to $\psi\star\mu(t)$ uniformly in $t\in\R^d$. Hence the last function is almost periodic too,
and $\mu$ is an almost periodic measure. \bs

{\bf Proof of Theorem \ref{T1}}. Let $\phi$ is a $C^\infty$-function  with compact support in $\R^d$ and non-negative Fourier transform $\Psi(x)$ (to construct such $\phi$
 should be taken $\phi(t)=\int\bar\p(s-t)\p(s)ds$ for $C^\infty$-function $\p$ with compact support).  There exists a ball $B(a,\rho)$ such that $\inf_{x\in B(a,\rho) }\Psi(x)\ge\eta>0$. Therefore,
$$
\eta\mu(B(a,\rho))\le\int_{B(a,\rho)}\Psi(t)\mu(dt)\le\int_{\R^d}\Psi(t)\mu(dt)=\left|\int_{\R^d}\phi(y)\hat\mu(dy)\right|\le\int_{\supp\phi}|\phi(y)||\hat\mu|(dy).
  $$
  Then for any $x\in\R^d$
  $$
\eta\mu(B(x,\rho))\le\int_{\R^d}\Psi(t+a-x)\mu(dt)=\left|\int_{\R^d}\phi(y)e^{2\pi iy(a-x)}\hat\mu(dy)\right|\le\int_{\supp\phi}|\phi(y)||\hat\mu|(dy).
 $$
 Hence, values $\mu(B(x,\rho)))$ are uniformly bounded, and the measure $\mu$ is translation bounded. Proposition \ref{P4} implies  almost periodicity of $\mu$.
Since $\mu$ is translation bounded, we see that $\mu$ is a temperate distribution, hence by Proposition \ref{P3},  the convolution $\hat\mu\star\p$ is almost periodic for every $\p\in S(\R^d)$.  \bs
\medskip

{\bf Proof of Theorem \ref{T2}}. Assume that $\mu\not\equiv\nu$. Then there is a point $x\in\R^d$ and a continuous function $\p,\,0\le\p\le1$, with $\supp\p\subset B(0,1)$ such that $\mu\star\p(x)\neq\nu\star\p(x)$.
Let $2\e<|\mu\star\p(x)-\nu\star\p(x)|$. The function $g=\mu\star\p-\nu\star\p$ is almost periodic, hence for $r_\e$ large enough every ball $B(x',r_\e)$  contains the ball $B(x+\tau,1)$
 for some $\e$-almost period $\tau$ of the function $g$.
Therefore, $|g(x+\tau)|>\e$ and
$$
|\mu-\nu|B(x',r_\e)\ge|\mu-\nu|\star\p(x+\tau)
\ge|\mu\star\p(x+\tau)-\nu\star\p(x+\tau)|\ge\e.
$$
Note that every ball of radius $R$ contains at least $\k_d(R/r_\e)^d$
disjointed balls of radius $r_\e$, where the constant $\k_d$ depends only on the dimension $d$. Therefore,
$$
|\mu-\nu|B(x_n,R_n)\ge \k_dr_\e^{-d}\e  R_n^d.
$$
We get a contradiction with conditions of the theorem.

The same arguments work in the case $\hat\mu\not\equiv\hat\nu$ after replacing a continuous $\p$ with $C^\infty$-function $\p$. \bs

\medskip

{\bf Proof of Theorem \ref{T3}}. Assume $\mu\not\equiv\nu$. Let $x$, $\p$, $\e$, $r_\e$ be the same as above. Take appropriate $\tau$ such that $B(x+\tau,r_\e)\subset B(x_n,R_n)$ for some $n$ and
\begin{equation}\label{dif}
|\mu\star\p(x+\tau)-\nu\star\p(x+\tau)|\ge\e.
\end{equation}
Let $M_\tau$ be a set of pairs $\{(a_n,c_n)\}$ such that either $a_n\in\supp\mu\cap B(x+\tau,1)$ or $c_n\in\supp\nu\cap B(x+\tau,1)$. By Theorem \ref{T1}, the measures $\mu$ and $\nu$ are translation bounded.
Therefore  quantities $\# M_\tau$ are bounded by a constant, which does not depend on $\tau$. We have
$$
|\mu\star\p(x+\tau)-\nu\star\p(x+\tau)|\le\sum_{(a_n,c_n)\in M_\tau}
|\p(a_n)-\p(c_n)|\le(\# M_\tau)\k_d N_1(\p)\max_{(a_n,c_n)\in M_\tau}|a_n-c_n|.
$$
It follows from conditions of the theorem that this inequality contradicts to \eqref{dif} for $n$ large enough. \bs
 \medskip

{\bf Proof of Theorem \ref{T4}}.  Assume that there is a point $x\in\R^d$ such that $0<\mu(x)<c$. Set $\e<(1/2)\min\{(c-\mu(x)),\mu(x)\}$ and $r<b/2$
such that $\mu(B(x,2r))<\mu(x)+\e$. Let $\p(x),\ 0\le\p(x)\le1$ be a continuous function with support in $B(0,r)$ such that $\p(0)=1$. We have
$$
\mu(x)\le\mu\star\p(x)\le\mu(B(x,r))<\mu(x)+\e.
$$
  For $R_n$ large enough the ball $B(x_n,R_n)\subset E$ contains the ball $B(x+\tau,r)$
for some $\e$-almost period $\tau$ of the function $\mu\star\p$.  Then $\mu\star\p(x+\tau)>0$, hence the ball $B(x+\tau,r)$ contains exactly  one point  $x'\in\supp\mu$ such that
$c\le\mu(x')\le\mu\star\p(x')$. Therefore,
$$
 c-\e<\mu\star\p(x'-\tau)\le\mu(B(x'-\tau,r))\le\mu(B(x,2r))<\mu(x)+\e,
$$
which contradicts our choice of $\e$. Consequently, the estimate $\mu(x)\ge c$ holds for all $x\in\supp\mu$.

Assume that there are points $x,\,x+h\in\supp\mu$ such that $|h|<b$. Let $\p,\ 0\le\p(x)\le1$ be a continuous function with support in $B(0,\e)$ with $\e<\min\{(b-|h|),c\}$ such that $\p(0)=1$. Then
$\mu\star\p(x)\ge c$ and $\mu\star\p(x+h)\ge c$.  For $R_n$ large enough  the ball $B(x_n,R_n)\subset E$
contains the ball $B(x+\tau,|h|+\e)$ for some $\e$-almost period $\tau$ of the function $\mu\star\p$. Then
$$
\mu\star\p(x+\tau)>0,\qquad\mu\star\p(x+h+\tau)>0,
$$
hence there are points $x'\in B(x+\tau,\e),\,x''\in B(x+h+\tau,\e)$ such that $\mu(x')>0,\,\mu(x'')>0$. Since $|x'-x''|<|h|+2\e<b$, we get a contradiction.
\smallskip

In the case of a complex measure $\mu$ condition \eqref{mu} and Proposition \ref{P3} yield that $\mu$ is an almost periodic distribution. From here, repeating the previous arguments with replacing
$\mu$ by $|\mu|$, we obtain that inequalities $|\mu(\l)|\ge c$ and $|\l-\l'|\ge\b>0,\,\l\neq\l'$ hold for all points $\l,\,\l'\in\supp\mu$. Under these conditions
 the first equality in \eqref{m}  was proved in \cite[Theorem 5]{F4}, and the second one is contained in \cite[Theorem 4]{F3}.   \bs

\bigskip
 \section{Proof of  sufficiency in  Theorem \ref{T5}}\label{S4}
\bigskip

Let $A=\{a_n\}\subset\R$ be a sequence such that the Fourier transform $\hat\mu_A$ is measure \eqref{b}. Our theorem does not depend on shifts along real axis, hence we can assume
 $a_n\neq0$ for all $n$.  Since $\mu_A$ is a non-negative Poisson measure, we see that $A$ is an almost periodic set.  By Proposition \ref{P2}, we can renumber $a_n$ such that $a_n=n/D+\phi(n),\,n\in\Z$
 with a bounded function $\phi(n)$ and a density $D$, and we get condition \eqref{l1}.

Consider the sum
\begin{equation}\label{sum}
\sum_{n\in\Z}\frac{1}{x+iy-a_n}=\sum_{n\in\Z}\frac{1}{x+iy-n/D-\phi(n)}.
\end{equation}
Set $M=\sup_\R|\phi(n)|$. Since
\begin{equation}\label{sum2}
\sum_{|n/D|<M+|x|+1}\left|\frac{1}{x+iy-n/D-\phi(n)}\right|+\sum_{|n/D|\ge M+|x|+1}\left|\frac{1}{x+iy-n/D-\phi(n)}+\frac{1}{x+iy+n/D-\phi(-n)|}\right|\le R,
\end{equation}
with $R$ depending on $x$ and $y$, we get that \eqref{sum} is finite for every $x\in\R,\,y>0$.

Furthermore, the sums
$$
 \sum_{n\in\Z,|a_n|<N}\frac{1}{a_n} \quad\mbox{and}\quad \sum_{n\in\Z,|n/D|<N}\frac{1}{a_n}=\frac{1}{a_0}+\sum_{n=1}^N\left[\frac{1}{n/D+\phi(n)}+\frac{1}{-n/D+\phi(-n)}\right]
$$
 differ for a uniformly bounded (with respect to $N$) number of terms, and each of these terms tends to $0$ as $N\to\infty$. Therefore these sums have the same finite limit $L$ as $N\to\infty$,
and we get \eqref{l2}. By Lindelof's Theorem, the function $F(z)$ from \eqref{F} has exponential growth. Moreover, we obtain that the sum
$$
   \sum_{n\in\N}\left|\frac{1}{a_n}+\frac{1}{a_{-n}}\right|
$$
is also converges, hence
$$
  F(z)=(1-z/a_0)e^{z/a_0}\prod_{n\in\N}(1-z/a_n)(1-z/a_{-n})e^{z(1/a_n-1/a_{-n})}=e^{Lz}f(z),
$$
where
\begin{equation}\label{f}
  f(z)=(1-z/a_0)\prod_{n\in\N} (1-z/a_n)(1-z/a_{-n})
\end{equation}
is a well-defined  entire function of exponential growth.
\medskip

Set for $\Im z>0$
$$
\xi_z(t)=\begin{cases}-2\pi ie^{2\pi itz} &\text{if }t>0,\\0&\text{if }t\le0,\end{cases}
$$
It is not hard to check that in the sense of distributions $\hat\xi_z(\l)=1/(z-\l)$.
Let $\p(t)$ be any  even nonnegative $C^\infty$-function such that
$$
\supp\p\subset(-1,1)\quad\text{and}\quad \int\p(t)dt=1.
$$
 Set $\p_\e(t)=\e^{-1}\p(t/\e)$ for $\e>0$.  The functions $\xi_z(t)\star\p_\e(t)$ and $\hat\xi_z(\l)\hat\p_\e(\l)$ belong to $S(\R)$. Therefore,
\begin{equation*}
   (\mu_A,\hat\xi_z(\l)\hat\p_\e(\l))=(\hat\mu_A,\xi_z\star\p_\e(t)).
\end{equation*}
We have
\begin{equation*}
 (\mu_A(\l),\hat\xi_z(\l)\hat\p_\e(\l))=\frac{\hat\p(\e a_0)}{z-a_0}+\sum_{n\in\N}\left[\frac{\hat\p(\e a_n)}{z-a_n}+\frac{\hat\p(\e a_{-n})}{z-a_{-n}}\right].
\end{equation*}
Clearly, $\hat\p_\e(t)\to1$ as $\e\to0$. Taking into account \eqref{sum2}, we can apply Lebesgue's Dominate Convergence Theorem. Hence,
\begin{equation}\label{q1}
(\mu_A,\hat\xi_z(\l)\hat\p_\e(\l))\to \frac{1}{z-a_0}+\sum_{n\in\N}\left[\frac{1}{z-a_n}+\frac{1}{z-a_{-n}}\right]\qquad(\e\to0).
\end{equation}

Set $n(s)=\sum_{\g\in\G:\,0<\g<s}|b_\g|=|\hat\mu_A|(0,s)$. If $|\hat\mu_A|(0,s)\le Ce^{2\pi sK}$, we obtain for $y>K$
\begin{equation}\label{i1}
   \sum_{\g\ge r}|b_\g|e^{-2\pi\g y}=\int_r^\infty e^{-2\pi sy}n(ds)\le\lim_{T\to\infty}n(T)e^{-2\pi Ty}+2\pi y\int_r^\infty e^{-2\pi sy}n(s)ds<\infty.
\end{equation}
Therefore the series $\sum_{\g>0}b_\g e^{2\pi i\g z}$ absolutely and uniformly  converges for every fixed $y=\Im z>K$.

Further,  we have for $\e>0$
\begin{multline*}
 \frac{i}{2\pi}(\hat\mu_A(t),\xi_z\star\p_\e(t))= \sum_{\g<0}b_\g e^{2\pi i\g z}\int_{s<\g} e^{-2\pi isz}\p_\e(s)ds+b_0\int_{s<0}e^{-2\pi isz}\p_\e(s)ds\\
+\sum_{\g>0}b_\g e^{2\pi i\g z}\int_{s<\g} e^{-2\pi isz}\p_\e(s)ds.
 \end{multline*}
It is easy to see that all integrals  are bounded and
$$
\int_{s<\g<0} e^{-2\pi isz}\p_\e(s)ds\to0,\quad \int_{s<0}e^{-2\pi isz}\p_\e(s)ds\to\frac{1}{2}, \quad \int_{0<s<\g}e^{-2\pi isz}\p_\e(s)ds\to1\qquad(\e\to0).
$$
Using Lebesgue's Dominate Convergence Theorem, we obtain
$$
  (\hat\mu_A,\xi_z\star\p_\e(t))\to-2\pi i\sum_{\g\in\G\cap(0,+\infty)}b_\g e^{2\pi i\g z}-\pi ib_0\qquad(\e\to0),
$$
and by \eqref{q1},
\begin{equation}\label{der}
f'(z)/f(z)=\frac{1}{z-a_0}+\sum_{n\in\N}\left[\frac{1}{z-a_n}+\frac{1}{z-a_{-n}}\right]=-2\pi i\sum_{\g\in\G\cap(0,+\infty)}b_\g e^{2\pi i\g z}-\pi iD,
\end{equation}
where we replace $b_0$ with the density $D$. Therefore, $f'(z)/f(z)$ is an absolutely convergent Dirichlet series on the line $x\in\R,\,y=\const>K$, and its mean value equals $-\pi iD$.

Changing the order of summation and integration in \eqref{der},  we get
\begin{equation}\label{pr}
  \log f(x+iy)-\log f(iy)=\int_{iy}^{x+iy}\frac{f'(\zeta)}{f(\zeta)}d\zeta=-\sum_{\g\in\G\cap(0,+\infty)}b_\g \frac{(e^{2\pi i\g x}-1)e^{-2\pi\g y}}{\g}-iD\pi x.
\end{equation}
Since $|e^{2\pi i\g x}-1|/|\g|\le 2\pi |x|$, we see that \eqref{pr} is well defined for all $x\in\R$ and $y>K$. Then
$$
\|\log f(x+iy)+iD\pi x\|_W\le \sum_{0<\g<1}\frac{|b_\g|}{\g}2e^{-2\pi\g y}+\sum_{\g\ge1}\frac{|b_\g|}{\g}2e^{-2\pi\g y}+|\log|f(iy)||<\infty.
$$
Condition \eqref{int} implies that the first sum in the right-hand side converges, and \eqref{i1} implies that the second one converges as well.

 Redefine  $y=y_0$. We obtain
$$
   f(x+iy_0)e^{iD\pi x}=\sum_{\o\in\O}\b_\o e^{2\pi i\o x},\quad \sum_{\o\in\O}|\b_\o|=\|f(x+iy_0)e^{iD\pi x}\|_W\le e^{\|\log f(x+iy_0)+iD\pi x\|_W}<\infty,
 $$
with $\b_\o\in\C$ and a countable spectrum $\O$.  The entire function $f(z+iy_0)e^{iD\pi z}$ has the exponential growth with zeros in a horizontal strip of bounded width.
It follows from \cite[\S 1, Ch.VI]{L}, that $\O$ is bounded. Hence the function
$$
    f(z)=\sum_{\o\in\O}\b_\o e^{\pi(2\o-D)y_0}e^{\pi i(2\o-D)z}
$$
is also a Dirichlet series and $\|f\|_W<\infty$.  \bs

\bigskip
 \section{Proof of  necessity in Theorem \ref{T5}}\label{S5}
\bigskip

Let $Q$ be  Dirichlet series of the form \eqref{S}
 with the zero set $A\subset\R$. Since zeros of $Q(z)$ lie in a strip of finite width, we get that $\k:=\inf\O\in\O$ and $\k':=\sup\O\in\O$
 (see \cite{L}, Ch.VI, Cor.2).  We obtain for $z=x+iy$
$$
   \log Q(x+iy)=\log q_\k+2\pi i\k(x+iy)+\log (1+H(x+iy)),\quad H(x+iy)=\sum_{\o\in\O\setminus\{\k\}}\frac{q_\o}{q_\k e^{2\pi(\o-\k)y}}e^{2\pi i(\o-\k)x}.
$$
Taking into account that $\sum_\o|q_\o|<\infty$, we can take a finite number of points  $\o_1,\dots,\o_N\in\O$ and then $s>0$ such that
$$
\sum_{\o\in\O\setminus\{\k,\o_1,\dots,\o_N\}}|q_\o/q_\k|<1/3,\qquad \sum_{j=1}^N e^{-2\pi(\o_j-\k)s} |q_{\o_j}/q_\k|<1/3.
$$
Also, we may suppose that for the same $\o_1,\dots,\o_N$ and  $s$
\begin{equation}\label{n}
\sum_{\o\in\O\setminus\{\k',\o_1,\dots,\o_N\}}|q_\o/q_{\k'}|<1/3,\qquad \sum_{j=1}^N e^{2\pi(\o_j-\k')s} |q_{\o_j}/q_{\k'}|<1/3.
\end{equation}

So $\|H(x+iy)\|_W<2/3$ for $y\ge s$, and the function
$$
\log(1+H(x+iy))=\sum_{n=1}^\infty\frac{(-1)^{n+1}H^n(x+iy)}{n}
$$
 is well-defined and belongs to $\H$.  Since $\spec H(x+iy)\subset(0,+\infty)$, we see that spectrum of the function
 \begin{equation}\label{logQ}
  \log Q(x+iy)-2\pi i\k(x+iy)-\log q_\k=\log(1+ H(x+iy))
\end{equation}
lies in $(0,+\infty)$ too.

From  almost periodicity of $Q$ it follows that for any $\e>0,\,s>0$ there exists $m(\e,s)>0$ such that
 $$
 |Q(z)|\ge m(\e,s)\quad\mbox{for}\quad |\Im z|\le s\quad\mbox{and}\quad z\not\in A(\e):=\{z:\,\dist(z,A)<\e\}.
 $$
 (\cite{L}, Ch.6, Lemma 1). Hence we have for any fixed $y>0$
 \begin{equation}\label{be}
\inf_{x\in\R}|Q(x\pm iy)|>0.
\end{equation}
 Further, the zero set $A$ of the Dirichlet series $Q(z)$ is almost periodic, and Proposition \ref{P1} implies that each connected component of $A(\e)$  contains
 no segment of length $1$ for sufficiently  small $\e$.
 Hence there are sequences $L_k\to+\infty,\,L'_k\to-\infty$ and $m>0$ such that
$$
   |Q(x+iy)|>m\quad\text{for}\quad x=L_k \quad\text{or}\quad x=L'_k,\quad |y|\le s.
$$
Let $\p$ be $C^\infty$-function with compact support on some interval $(-r,r)$. Set
$$
\Phi(z)=\int_{-\infty}^\infty \p(t)e^{-2\pi itz}dt.
$$
Clearly,  $\Phi(x+iy)$ is an entire function, which  equals the Fourier transform of the function $\p(t)e^{2\pi ty}$ for any fixed $y$.
Therefore, $\Phi(x+iy)$  belongs to $S(\R)$ for each fixed $y\in\R$, and  we have for its converse Fourier transform
\begin{equation}\label{cF}
 \check\Phi(t+iy)=\p(t)e^{2\pi ty},\qquad t\in\R.
\end{equation}
Hence the function $\Phi(x+iy)$  tends to zero for $x=L_k$ or $x=L'_k$ as $k\to\infty$ uniformly with respect to $y\in[-s,s]$.
 Consider  integrals of the function $\Phi(z)Q'(z)Q^{-1}(z)$ over the boundaries of  rectangles
$\Pi_k=\{z=x+iy:\,L'_k<x<L_k,|y|<s\}$.  These integrals  tend to the difference
\begin{equation}\label{in}
   \int_{-is-\infty}^{-is+\infty}\Phi(z)Q'(z)Q^{-1}(z)dz-\int_{is-\infty}^{is+\infty}\Phi(z)Q'(z)Q^{-1}(z)dz=:I_1-I_2
\end{equation}
 as $k\to\infty$. On the other hand, they tend to the sum
  \begin{equation}\label{r}
  2\pi i\sum_{\l:Q(\l)=0}\Res_\l \Phi(z)Q'(z)Q^{-1}(z)=2\pi i\sum_{\l:Q(\l)=0}a(\l)\hat\p(\l)=2\pi i(\mu_A,\hat\p)=2\pi i(\hat\mu_A,\p),
\end{equation}
where $a(\l)$ is the multiplicity of zero of $Q(z)$ at the point $\l$.

We have $Q'/Q=(\log Q)'$, hence \eqref{logQ} implies
 \begin{equation}\label{QQ}
 Q'(x+is)/Q(x+is)=\sum_{\g\in\G^*} p_\g e^{2\pi i\g x},\quad \sum_{\g\in\G^*}|p_\g|<\infty,\quad p_0=2\pi i\k,\quad p_\g=0\quad\text{for}\quad\g<0.
 \end{equation}
 It follows from \eqref{n} that the same arguments show
$$
 Q'(x-is)/Q(x-is)=\sum_{\g\in\G_*} \tilde p_\g e^{2\pi i\g x},\quad \sum_{\g\in\G_*}|\tilde p_\g|<\infty,\quad\tilde p_0=2\pi i\k',\quad\tilde p_\g=0\quad\text{for}\quad\g>0.
$$
   The functions $\Phi(x\pm is)$ belong to $S(\R)$, hence $|x|^2\Phi(x\pm is)\to0$ as $|x|\to\infty$.
 Changing the order of integration and summation and taking into account \eqref{cF}, we obtain for the difference of integrals  \eqref{in}
\begin{multline}\label{i}
 I_1-I_2=\sum_{\g\in\G_*}\tilde p_\g\int_{-\infty}^{+\infty}\Phi(x-is)e^{2\pi i\g x}dx-\sum_{\g\in\G^*}p_\g\int_{-\infty}^{+\infty}\Phi(x+is)e^{2\pi i\g x}dx\\
=\sum_{\g\in\G_*}\tilde p_\g e^{-2\pi\g s}\p(\g)-\sum_{\g\in\G^*}p_\g e^{2\pi\g s}\p(\g).
\end{multline}
Since $\supp\p\subset(-r,r)$, we obtain
\begin{equation}\label{b0}
|(\hat\mu_A,\p)|\le (2\pi)^{-1}\left(\sum_{\g\in\G^*}|p_\g|+\sum_{\g\in\G_*}|\tilde p_\g|\right)e^{2\pi rs}\sup_{|y|\le r}|\p(y)|.
\end{equation}

Hence the distribution $\hat\mu_A$ has a unique expansion to a linear functional on the space of continuous functions $g$ on $[-r,r]$ such that $g(-r)=g(r)=0$  with the same bound.
Since we can expand this functional to the space of all continuous functions on $[-r,r]$ with the same bound  too, we see that $\hat\mu_A$ is a complex measure.
 Then \eqref{r} and \eqref{i} implies
 \begin{equation}\label{hatmu}
\hat\mu_A=\sum_{\g\in\G_*}\frac{\tilde p_\g e^{-2\pi\g s}}{2\pi i}\d_\g-\sum_{\g\in\G^*}\frac{p_\g e^{2\pi\g s}}{2\pi i}\d_\g.
\end{equation}
Therefore, $\hat\mu_A$ is a purely point measure. Also, it follows from \eqref{b0} that
\begin{equation}\label{h}
|\hat\mu_A|(-r,r)\le C(s)e^{2\pi sr}.
\end{equation}
 Furthermore, the functions $Q(z)$ and $f(z)$ from \eqref{f} have exponential growth and the same zeros. Therefore,
 $$
Q(z)=Ce^{\a z+i\b z}f(z),\qquad C\in\C,\quad \a,\,\b\in\R.
$$
Hence,
\begin{equation}\label{x}
  \frac{Q'(x+iy)}{Q(x+iy)}-\frac{f'(x+iy)}{f(x+iy)}=\a+i\b.
\end{equation}
It follows from \eqref{h}  that the function $f'(x+iy)/f(x+iy)$  holds \eqref{der} for $y>s$.  Let
$$
\log(1+H(x+iy))=\sum_{\g\in\G_1}c_\g e^{2\pi i\g x},\quad c_\g=c_\g(y),\quad\G_1\subset(0,+\infty),\quad\sum_{\g\in\G_1}|c_\g|<\infty.
$$
By \eqref{logQ},
$$
  Q'(x+iy)/Q(x+iy)=(\log Q(x+iy))'=\sum_{\g\in\G_1} 2\pi i\g c_\g e^{2\pi i\g x}+2\pi i\k.
$$
Therefore, \eqref{x} and \eqref{der} yield
$$
\a=0,\quad \b=\pi i(D+2\k),\quad 2\pi i\g c_\g=-2\pi i e^{-2\pi\g y}b_\g,\quad \G_1=\G\cap (0,+\infty),
$$
 and
$$
 \sum_{0<\g<1}|b_\g/\g|\le e^{2\pi y}\sum_{0<\g<1}|c_\g|<\infty.
$$
We obtain \eqref{int}.

By \eqref{be}, $1/Q(x\pm i\eta)\in\H$ for every $\eta>0$. The spectrum $\O$ is bounded, hence $Q'(x\pm i\eta)\in\H$ as well.
Therefore, $Q'(x\pm i\eta)/Q(x\pm i\eta)\in\H$, and we get the representation
$$
 Q'(x\pm i\eta)/Q(x\pm i\eta)=\sum_{\g\in S} h_\g^\pm e^{2\pi i\g x},\qquad \sum_{\g\in S}|h_\g^\pm|<\infty,
$$
where $S$ is a countable subsets of $\R$. Arguing as above, instead of \eqref{hatmu} we obtain the representation
$$
\hat\mu_A=\sum_{\g\in  S}\frac{h_\g^- e^{-2\pi\g\eta}-h_\g^+e^{2\pi\g\eta}}{2\pi i}\d_\g.
$$
Therefore, \eqref{h} is valid for all $s=\eta>0$. This gives \eqref{var}. \bs
\bigskip

I thank  Professor Szilard Revesz from Renyi Institute of Mathematics for the hospitality and useful discussions.


\begin{thebibliography}{19}


\bibitem{B} Bohr, H. Almost Periodic Functions, ed. Chelsea,
New-York, 1951.

\bibitem{D} Directions in Mathematical Quasicrystals, M.Baake, R.Moody, eds.
CRM Monograph series 2000 {\bf 13}, AMS, Providence RI, 379p.

\bibitem{F0} Favorov, S.Yu. Zeros of holomorphic almost periodic functions //Journal
d'Analyse Mathematique, v.84, (2001), p.51-66.

\bibitem{Fa} Favorov, S.Yu. Almost periodic divisors, holomorphic functions,
and holomorphic mappings// Bulletin des Sciences Mathematiques, v.127 (2003),
859-883.

\bibitem{F1} Favorov, S.Yu. Uniqueness Theorems for Fourier Quasicrystals and Temperate Distributions with Discrete Support. Proc. Amer. Math. Soc. 149 (2021), 4431-4440.

\bibitem{F2} Favorov, S.Yu. The crystalline measure that is not a Fourier quasicrystal, arXiv:2401.01121 (2024)

\bibitem{F3} Favorov, S.Yu. {\it Local Wiener's Theorem and Coherent Sets of Frequencies}. Analysis Math., {\bf 46} (4) (2020), 737–746
DOI: 10.1007/s10476-020-0042-x

\bibitem{F4} Favorov, S.Yu. {\it Large Fourier Quasicrystals and Wiener's Theorem},  Journal of Fourier Analysis and Applications, {\bf 25}, Issue 2, (2019), 377-392,
 DOI 10.1007/s00041-017-9576-0

\bibitem{F} Favorov, S.Yu. Generalized Fourier quasicrystals, almost periodic sets, and zeros of  Dirichlet series.
arXiv:2311.02728 (2023).

\bibitem{FK} Favorov, S.Yu. Kolbasina, Ye.Yu. Almost periodic discrete sets// Journal of
Mathematical Physics, Analysis, Geometry. v.6 (2010), No.1, 1-14.

\bibitem{FK1} Favorov, S.Yu., Ye.Yu. Kolbasina. Perturbations of discrete lattices and
almost periodic sets// Algebra and Discrete Mathematica, 2010, v.9, No.2, 48-58.

\bibitem{FRR} Favorov, S.Yu., Rashkovskii, A.Yu. and Ronkin, L.I. Almost
periodic divisors in a strip //Journal d'Analyse Mathematique, Vol.74 (1998), 325-345.

\bibitem{FRR1} Favorov, S.Yu., Rashkovskii, A.Yu. and Ronkin, L.I. Almost periodic currents and holomorphic chains //Comptes Rendus Acad.
Sci. Paris, t. 327 (1998), Serie I, p. 302-307.

\bibitem{K} P.Koosis, The logarithmic integral, Vol.I. Cambridge university press, Cambridge,  New York, New Rochelle, Melburn, Sydney, 1988.

\bibitem{KS} P.Kurasov, P.Sarnak, Stable polynomials and crystalline measures. J. Math. Phys. 61, no. 8. 083501 (2020); https://doi.org/10.1063/5.0012286)

\bibitem{L} Levin, B.Ja. Distributions of Zeros of Entire Functions. Transl. of Math. Monograph, Vol.5, AMS Providence, R1, 1980.

\bibitem{Lag1} Lagarias, J.C. {\it Mathematical Quasicrystals and the Problem of Diffraction}, in \cite{D}, 61-93.

\bibitem{Law} Lawton, W.M. Tsikh, A.K. {\it Fourier Quasicrystals on $\R^n$}

\bibitem{M} Meyer, Y. {\it Quasicrystals, Almost Periodic Patterns, Mean--periodic Functions, and Irregular Sampling},
Afr. Diaspora J. Math., {\bf 13} (1) (2012), 1--45.

\bibitem{M1} Meyer, Y. {\it  Global and local estimates
on trigonometric sums}, Trans. R. Norw. Soc. Sci. Lett. 2018(2) 1-25.

\bibitem{M2} Meyer, Y. {\it Multidimensional crystalline measures} Trans. R. Norw. Soc. Sci. Lett. 2023(1) 1–24.

\bibitem{OU} Olevskii, A., Ulanovskii A.M.,  Fourier quasicrystals with unit masses// Comptes Rendus Mathématique, 2020, 358, no 11-12, p. 1207-1211
https://doi.org/10.5802/crmath.142

\bibitem{OU1} Olevskii, A., Ulanovskii A.M., A Simple Crystalline Measure.
arXiv:2006.12037v2, (2020).

\bibitem{Q} Quasicrystals and Discrete Geometry. J.Patera,ed., Fields Institute Monographs 1998, AMS, Providence RI, 289p.

\bibitem{R} Ronkin, L.I. {\it Almost Periodic Distributions and Divisors in Tube
Domains,} Zap. Nauchn. Sem. POMI {\bf 247} (1997)  210--236 (Russian).

\bibitem{R1} Rudin, W.: Fourier Analysis on Groups. Interscience Publications, a Division of John Wiley and Sons, New York  (1962)

\bibitem{T} Tornehave, H. {\it On the zeros of entire almost periodic function}.
The Harald Bohr Centenary (Copenhagen 1987). Math. Fys. Medd. Danske, {\bf 42}, no.3 (1989), 125-142.


\end{thebibliography}
\end{document}